\documentclass{article}
\usepackage[english]{babel}
\usepackage[cp1251]{inputenc}
\usepackage{amssymb,amsmath,amsthm}
\newtheorem{thm}{Theorem}
\newtheorem*{thm*}{Theorem}
\newtheorem*{lem*}{Lemma}
\newtheorem{lem}{Lemma}
\newtheorem{cor}{Corollary}
\newtheorem{pr}{Proposition}

\begin{document}
\title{The necessary and sufficient condition for solvability of a partial integral equation}
\author{\sc{Eshkabilov Yu.Kh.}\\
\it{National University of Uzbekistan}\\
e-mail: yusup62@rambler.ru}
\date{}
\maketitle
\begin{abstract}

Let $T_1: L_2(\Omega^2) \to L_2(\Omega^2)$ be a partial integral
operator with the kernel from $C(\Omega^3)$ where $\Omega=[a,b
]^\nu.$ In this paper we investigate solvability of a partial
integral equation $f-\varkappa T_1 f=g_0$ in the space
$L_2(\Omega^2)$ in the case when $\varkappa$ is a characteristic
number. We proved the theorem describing the necessary and
sufficient condition for solvability of the partial integral
equation $f-\varkappa T_1 f=g_0.$


\end{abstract}

In the models of solid state physics [1] and also in the lattice
field theory [2], there appear so called discrete Schrodinger
operators, which are lattice analogues of usual Schrodinger
operators in a continuous space. The study of spectra of lattice
Hamiltonians (that is discrete Schrodinger operators) is an
important matter of mathematical physics. Nevertheless, on
studying spectral properties of discrete Schrodinger operators
three appear partial integral  equations in a Hilbert space of
multi- variable functions [1,3]. Therefore, on the investigation
of spectra of Hamiltionians considered on a lattice, the study of
a solvability problem for partial integral equations in  $L_2$ is
essential (and even interesting from the point of view of
functional analysis).

A question on the existence of a solution of partial integral
equation (PIE) for functions of two variables was considered in
[4-8] and other works. In the work by author [9], the PIE
$f-\varkappa T_1 f=g_0$ was studied in the space $L_2(\Omega^2),$
where $\Omega=[a,b]^\nu,$ for a partial integral operator (PIO)
$T_1: L_2(\Omega^2)\to L_2(\Omega^2)$ with the kernel $k(x,s,y)$
being a three variables continuous function on $\Omega^3.$ The
concept of determinant for the PIE as a continuous function on
$\Omega$ and concepts of regular number, singular number,
characteristic number and essential number for a PIE are given.
Theorems on solvability of the PIE are proved in the case when
$\varkappa$ is the regular and essential number [9]. In this paper
we study solvability of the PIE $f-\varkappa T_1 f=g_0$ when
$\varkappa$ is the characteristic number, i.e. the paper continues
the work by author [9].

Let $L^0=L^0(\Omega)$ be a space of classes of complex measurable
functions $b=b(y)$ on $\Omega$. We define in the space
$L^0(\Omega)$ the norm of an element $b\in L^0(\Omega)$ by the
equality $\| b\|=\sqrt{\int |b(t)|^2dt}.$

We denote by $L_{2,0}(\Omega^2)$ the totality of classes of
complex measurable functions $f(x,y)$ on $\Omega\times \Omega$
satisfying the condition: $\int |f(x,y)|^2dx$ exists for almost
all $y\in\Omega.$ It is easy to note that $L_{2,0}(\Omega^2)$ is a
linear space over $\mathbb C$ and $L_2(\Omega^2) \subset
L_{2,0}(\Omega^2).$ For each $b(y)\in L^0$ and $f(x,y)\in
L_{2,0}(\Omega^2),$ we define the function $b\circ f$ by the
formula $(b\circ f)(x,y)=b(y)f(x,y).$ Then for any $b\in L^0$ we
have $b\circ f\in L_{2,0}(\Omega^2),$ where $f\in
L_{2,0}(\Omega^2).$ For any $f,g\in L_{2,0}(\Omega^2),$ the
integral $\int f(x,t) g(x,t)dx$ exists at almost all $t\in\Omega$
and $\varphi(t)=\int f(x,t) g(x,t)dx\in L^0.$

Let $\nabla$ be the Boolean algebra of idempotents in $L^0.$ A
system $\{ f_1, f_2, \ldots, f_n\}\subset L_{2,0}(\Omega^2)$ is
called $\nabla$-{\it linearly independent,} if for all
$\pi\in\nabla$ and $b_1(y), b_2(y), \ldots, b_n(y)\in L^0$ from
$\sum\limits_{k=1}^n \pi\circ (b_k\circ f_k)=\theta$ it follows
$\pi\cdot b_1= \pi\cdot b_2= \cdots=\pi\cdot b_n=\theta$ [10,11].

Consider the mapping $\langle\cdot,\cdot\rangle:
L_{2,0}(\Omega^2)\times L_{2,0}(\Omega^2) \to L^0$ acting by the
rule
$$
\langle f,g\rangle = \int f(s,y)\overline {g(s,y)}ds,\quad f,g\in
L_{2,0}(\Omega^2).
$$

For every $b\in L^0,$ we have $\langle b\circ f,g\rangle = b\cdot
\langle f,g\rangle,$ where $f,g\in L_{2,0} (\Omega^2),$ i.e. the
mapping $\langle\cdot,\cdot\rangle$ satisfies the condition of
$L^0$-valued internal product [11].

In the space $\cal {H}$ $=L_{2,0}(\Omega^2)$, we consider a
partial integral operator (PIO) $S$ defined by
$$
Sf=\int\limits_\Omega q(x,s,y)f(s,y)ds,\quad f\in{\cal H}
$$
where $q(x,s,y)\in L_2(\Omega^3).$ The function $q(x,s,y)$ is
called \emph{kernel} of the PIO $S_1.$

A kernel $\overline{q(s,x,y)}$ corresponds to the adjoint operator
$S_1^*,$ i.e.
$$
S^*f=\int\limits_\Omega\overline{q(s,x,y)}f(s,y)ds,\quad f\in \cal
H.
$$

Consider a family of operators $\{ S_\alpha\}_{\alpha\in \Omega}$
in $L_2(\Omega)$ associated with $S_1$ by the following formula
$$
S_\alpha\varphi=\int\limits_\Omega q(x,s,\alpha)\varphi(s)ds,\quad
\varphi\in L_2(\Omega),
$$
where $q(x,s,y)$ is the kernel of $S.$

Further, if the set of integrability of an integral is absent, we
mean integrability by the set $\Omega.$

Now we consider the equation
\begin{equation}
f-\varkappa Sf=g_0,
\end{equation}
in the space $\cal H$ where $f$ is an unknown function from $\cal
H$, $g_0\in\cal H$ is a given function, $\varkappa\in\mathbb C$ is
a parameter of the equation.

For each $n\in\mathbb N,$ we define the measurable function
$$
\Pi^{(n)}= \Pi^{(n)} (x_1,\dots, x_n, s_1, \dots, s_n,\alpha)
$$
on $\Omega^n\times \Omega^n \times \Omega$ by means of an
$n$-ordinary determinant by the equality
$$
\Pi^{(n)} (x_1,\dots, x_n, s_1, \dots, s_n,\alpha)= \left|
\begin{array}{ccc} q(x_1,s_1,\alpha) & \dots & q(x_1, s_n, \alpha)
\\ \vdots & \vdots & \vdots \\ q(x_n,s_1,\alpha) & \dots & q(x_n, s_n,
\alpha) \\ \end{array} \right| .
$$

Now, for every $\varkappa\in\mathbb C$ we "formally"\ define
functions $D_1(y) =D_1(y;\varkappa)$ on $\Omega$ and $M_1(x,s,y) =
M_1(x,s,y;\varkappa)$ on $\Omega^3$ by means of the sum of
measurable functional series composed from sequences of measurable
functions $d_n(y)$ on $\Omega$ and $q_n(x,s,y)$ on $\Omega^3,$
respectively, by the following rules
$$
D_1(\alpha) = 1 + \sum\limits_{n\in\mathbb N} \frac{(-\varkappa)^
n}{n!} d_n(\alpha), \quad \alpha\in \Omega \eqno (a)
$$
and
$$
M_1(x,s,\alpha)= q(x,s,\alpha)\sum\limits_{n\in\mathbb N}
\frac{(-\varkappa)^ n}{n!} q_n(x,s,\alpha), \quad (x,s,\alpha) \in
\Omega^3 \eqno (b)
$$
where
$$
d_k(\alpha) = \int \dots \int \Pi^{(k)}(\xi_1,\dots, \xi_k, \xi_1,
\dots, \xi_k, \alpha)d\mu(\xi_1) \dots d\mu(\xi_k),
$$ $$
q_k(x,s,\alpha) = \int \dots \int \Pi^{(k+1)}(x,\xi_1,\dots,
\xi_k, s, \xi_1, \dots, \xi_k, \alpha)d\mu(\xi_1) \dots
d\mu(\xi_k).
$$

\begin{lem} For each $\varkappa \in \mathbb C$ the functions
$D_1(y)= D_1(y;\varkappa)$ $(a)$ and $M(x,s,y) = M(x,s,y;
\varkappa)$ $(b)$ are measurable on $\Omega$ and $\Omega^3,$
respectively. Moreover, for almost all $\alpha \in \Omega,$ there
exists the integral $\int \int |M_1(x,s,\alpha)|^2 dx ds.$
\end{lem}
\begin{proof}
It is clear that the operator $S$ is compact for almost all
$\alpha\in \Omega$. Let $\varkappa\in \mathbb C$ be an arbitrary
fixed number. We respectively denote as $\Delta_\alpha(
\varkappa)$ and $M_\alpha(x,s;\varkappa)$ the determinant and the
Fredholm minor of the operator $I-\varkappa S_\alpha$ for
$\alpha\in \Omega$ where $I$ is the identity operator in
$L_2(\Omega).$ Let $\varphi_n(y)$ and $\psi_n(x,s,y)$ be the
partial sums of the functional series $(a)$ and $(b),$
respectively. We have the sequences of measurable functions
$\varphi_n(y)$ on $\Omega$ and $\psi_n(x,s,y)$ on $\Omega^3$ such
that $\lim\limits_{n\to \infty} \varphi_n(y)= \Delta_y^{(1)}
(\varkappa) = D_(y;\varkappa)$ for almost all $y\in \Omega$ and
$\lim\limits_{n\to\infty} \psi_n(x,s,y)=M_y^{(1)}(x,s;\varkappa) =
M_1(x,s,y;\varkappa)$ for all $x,s \in \Omega$ and for almost all
$y\in \Omega.$ Therefore the function $D_1(y) = D_1(y;\varkappa)$
and the function $M_1(x,s,y)= M_1(x,s,y;\varkappa)$ are measurable
on $\Omega$ and $\Omega^3,$ respectively. It is known that if the
kernel $h(x,s)$ of the integral operator $A \varphi = \int h(x,s)
\varphi(s) ds,$ $\varphi\in L_2(\Omega)$ is an element of the
space $L_2(\Omega^2),$ then the minor $M(x,s;\varkappa)$ of the
operator $I -\varkappa A$ is also an element of the space $L_2(
\Omega^2)$ []. Hence we have
$$
\int\int |M_1(x,s,\alpha)|^2 dxds < \infty\quad \mbox{ for almost
all } \quad \alpha\in \Omega.
$$
\end{proof}

The measurable functions $D_1(y)= D_1(y;\varkappa)$ and $M_1(x,s,
y)= M_1(x,s,y;\varkappa)$ are respectively called the determinant
and the minor of the operator $E-\varkappa S,$ $\varkappa\in
\mathbb C,$ where $I$ is the identity operator in $L_{2,0}
(\Omega^2).$

\begin{lem}
Let $S:L_{2,0}(\Omega^2) \to L_{2,0} (\Omega^2)$ be a PIO with the
kernel $q\in L_2(\Omega^3).$ If the homogeneous equation $\varphi
-\varkappa S_\alpha \varphi = \theta,$ $\varkappa\in \mathbb C,$
has only the trivial solution in $L_2(\Omega)$ for almost all
$\alpha\in \Omega',$ then PIE (1) is solvable in the space
$L_{2,0}(\Omega^2)$ for every $g_0\in L_{2,0}(\Omega^2).$
\end{lem}

\begin{proof}
Let $\varkappa\in \mathbb C,$ $g_0(x,y)$ be an arbitrary function
from the space $L_{2,0}(\Omega^2).$ Let the homogeneous equation
$\varphi - \varkappa S_\alpha \varphi= \theta$ has only a trivial
solution in the space $L_2(\Omega)$ for almost all $\alpha \in
\Omega'$. Then $D_1(\alpha)= D_1(\alpha;\varkappa)\ne 0$ for
almost all $\alpha\in \Omega$ and the equation $\varphi(x) -
\varkappa (S_\alpha\varphi) (x) = h_\alpha(x)$ has a solution
$\varphi_\alpha(x) \in L_2(\Omega)$ for almost all $\alpha\in
\Omega'$ where $h_\alpha(x) = g_0(x,\alpha)\in L_2(\Omega).$
Moreover, the solution $\varphi_\alpha(x)$ has the form []
$$
\varphi_\alpha(x) = h_\alpha(x)  + \varkappa \int \frac{M_1(x,s,
\alpha;\varkappa)}{D_1(\alpha;\varkappa)} h_\alpha(s) ds.
$$
We have
$$
\int\int \left| \frac{M_1(x,s,\alpha;\varkappa)}{D_1(
\alpha;\varkappa)}\right| dx ds <\infty \ \mbox{ for almost all }
\ \alpha\in \Omega.
$$
It means that we can define a PIO $W=W(\varkappa):
L_{2,0}(\Omega^2) \to L_{2,0}(\Omega^2)$ with the kernel
$\frac{M_1(x,s,\alpha;\varkappa)}{D_1( \alpha;\varkappa)}$ [].
Therefore we have $f_0(x,y) = g_0(x,y) + \varkappa (Wg_0) (x,y)
\in L_{2,0}(\Omega^2)$ and $\varphi_\alpha(x) = f_0(x,\alpha)$ for
all $x\in \Omega$ and for almost all $\alpha\in \Omega.$ So the
function $f_0(x,y)$ is a solution of the equation (1).
\end{proof}

The following two propositions are proved analogously to
Propositions 1 and 2 from [].

\begin{pr}
Let $S:L_{2,0}(\Omega^2) \to L_{2,0}(\Omega^2)$ be a PIO with the
kernel $q\in L_2(\Omega^3).$ Then the following two conditions are
equivalent:

(i) a number $\lambda\in \mathbb C$ is an eigenvalue of the
operator $S;$

(ii) a number $\lambda\in \mathbb C$ is an eigenvalue of operators
$\{ S_\alpha\}_{\alpha\in \Omega_0}$ where $\Omega_0$ is a subset
of $\Omega$ such that $\mu(\Omega_0) >0.$
\end{pr}

\begin{pr}
If $\lambda\in \mathbb C$ is an eigenvalue of the PIO
$S:L_{2,0}(\Omega^2) \to L_{2,0}(\Omega^2)$ with the kernel
$q(x,s,y) \in L_2(\Omega^3),$ then the number $\overline \lambda$
is an eigenvalue of the operator $S^*.$
\end{pr}

\begin{thm}
Let $S:L_{2,0}(\Omega^2) \to L_{2,0}(\Omega^2)$ be a PIO with the
kernel $q(x,s,y)\in L_2(\Omega^3).$ Then every eigenvalue of the
PIO $S$ corresponds only to a finite number of $\nabla$-linearly
independent eigenfunctions.
\end{thm}
\begin{proof}
Let $\lambda\in\mathbb C$ be an eigenvalue of the PIO $S$ and
\begin{equation}
f_1,\ f_2,\ \ldots,\ f_m
\end{equation}
be some $\nabla$-linearly independent eigenfunctions, i.e.
\begin{equation}
\lambda f_j(x,y)=(Sf_j)(x,y),\quad j=1,2,\ldots,m.
\end{equation}
Since any linear combination of the eigenfunctions (2) of the
operator $S$ with coefficients from $L^0$ is also an
eigenfucntion, we can apply to the functions (2) the process of
$L^0$-orthogonalization [12]. Thus, we can assume that the
functions (4) are mutually orthogonal and normed in the sense of
$M^0$-valued internal products, i.e.
$$
\langle f_i,f_j\rangle =0,\quad i\ne j\,\ \mbox{ and }\ \langle
f_i,f_i\rangle =1.
$$
Therefore we can rewrite (3) in the following form
$$
\overline \lambda \cdot\overline {f_j(x,y)} = \left( S^* \overline
f_j\right) (x,y) = \int \overline{q(x,s,y)} \cdot
\overline{f_j(s,y)} ds.
$$
From here, it is easy to see that at almost all $x\in\Omega$ the
left side of this equality is a $L^0$-valued Fourier coefficient
of the function $\overline{q(x,s,y)}$ as and it is a function of
$(s,y)$ with respect to the orthogonal normed system (2). By the
Bessel inequality [12], one can write
$$
|\lambda|^2 \sum\limits_{j=1}^m |f_j(x,y)|^2 \le \int |q(s,x,y)|^2
ds\,\ \mbox{ for almost all }\ x\in \Omega.
$$
If we integrate both parts of this inequality by $x$ and $y$, we
obtain
$$
m\le |\lambda|^{-2} \int\int\int |q(x,s,y)|^2 dxdsdy <\infty.
$$
Hence, the number of $\nabla$-linearly independent functions,
corresponding to the eigenvalue $\lambda,$ is finite.
\end{proof}

Let $S$ be a PIO with the kernel $q(x,s,y)\in L_2(\Omega^2).$ A
number $\varkappa_0\in\mathbb C$ is called {\it a characteristic
value} of the PIE $f-\varkappa_0 Sf=g_0$ if the homogeneous
equation $f-\varkappa_0 Sf=0$ has a non-trivial solution. From
here, it is clear that any characteristic value $\varkappa_0$ of
the PIE $f-\varkappa Sf=g_0$ is non-zero.
\begin{cor}
Let $S:L_{2,0}(\Omega^2) \to L_{2,0}(\Omega^2)$ be a PIO with the
kernel $q(x,s,y)\in L_2(\Omega^3).$ Then any characteristic value
of the PIE $f-\varkappa Sf=g_0$ corresponds only to finite number
of $\nabla$-linearly independent eigenfunctions.
\end{cor}
\begin{thm}
Let $\varkappa$ be a characteristic numberof the PIE (1). Then the
homogeneous PIE
\begin{equation}
f-\varkappa Sf=0
\end{equation}
and the adjoint homogeneous PIE
$$
f-\overline{\varkappa} S^*f=0 \eqno (4')
$$
have the same number of $\nabla$-linearly independent solutions.
\end{thm}
\begin{proof}
Let $f_1,\ldots, f_m$ and $g_1,\ldots, g_n$ be $\nabla$-linearly
independent solutions of the the homogeneous equations (4) and
$(4')$, respectively. Assume that $m<n.$ We can suppose that
$f_1,\ldots, f_m$ and $g_1,\ldots, g_n$ are orthonormal systems in
the sense of $L^0$-valued internal product.

Define the function
$$
p(x,s,y)=q(x,s,y)-\sum\limits_{j=1}^m \overline{f_j(s,y)}g_j(x,y).
$$
We have $p(x,s,y)\in L_2(\Omega^3)$ since $f_j,g_k\in
L_2(\Omega^2).$ Consider two homogeneous PIE:
\begin{equation}
f-\varkappa Wf=0
\end{equation}
and
$$
f-\overline \varkappa W^*f=0 \eqno (5')
$$
where $W$ is the PIO with the kernel $p(x,s,y).$

Let $h(x,y)$ be a solution of the equation (5). Then we have $$
\langle h, g_j\rangle =\left\langle \varkappa Wh, g_j\right\rangle
= \left\langle h,\overline \varkappa S^*
g_j\right\rangle-\varkappa \left\langle h, f_j\right\rangle =
\left\langle h, g_j\right\rangle -\varkappa \left\langle h,
f_j\right\rangle,\ j=1,2,\ldots,m. $$ Hence, by virtue of
$\varkappa\ne 0,$
\begin{equation}
\left\langle h, f_j\right\rangle=0,\quad j=1,2,\ldots,m.
\end{equation}
Thus, any solution of the equation (5) satisfies the conditions
(6). But by virtue of this conditions, one can rewrite the
equation (5) in the form $f-\varkappa Sf=0,$ i.e. any solution of
the equation (5) satisfies the equation (4), too. We obtain that a
solution $h(x,y)$ of the eqaution (5) is to be in the form
$$
h(x,y)=\sum\limits_{j=1}^m \left( b_j\circ f_j\right) (x,y),\quad
b_j\in L^0,\quad j=1,2,\ldots, m.
$$
But we have $0=\langle h, f_k\rangle =\sum\limits_{j=1}^m
\left\langle b_j\circ f_j, f_k\right\rangle = \sum\limits_{j=1}^m
b_j\cdot \left\langle f_j, f_k\right\rangle = b_k,$ $k=1,2,\ldots,
m.$ Thus, we have $h(x,y)= \theta,$ i.e. the homogeneous PIE (5)
has only the trivial solution. We show that the adjoint equation
$(5')$ has non-trivial solutions. If we substitute
$g(x,y)=g_k(x,y),$ where $k>m,$ in the equation $(7')$ then we
obtain $g_k=\varkappa^* W^* g_k.$ Thus, we obtain the
contradiction to Proposition 2: the equation (5) has only the
trivial solution, but the adjoint equation $(5')$ has a
non-trivial solution. Hence the case $m<n$ is impossible. One can
prove similarly that the case $m>n$ is also impossible and we
obtain that $m=n.$
\end{proof}

\begin{thm} Let $\varkappa_0$ be a characteristic number of the PIE (1). Then:

a) the homogeneous equation $f-\varkappa_0Sf=0$ has a non-trivial
solution, moreover the set of all solutions of the homogeneous
equation is an infinite dimensional subspace of $\cal H;$

b) PIE (1) is solvable if and only if a given function $g_0$
satisfies the condition
$$
\langle g_0,g\rangle=0, \eqno (I)
$$
where $g\in\cal H$ is an arbitrary solution of the adjoint
homogeneous equation $f-\overline \varkappa_0S^*f=0.$
\end{thm}
\begin{proof} The proof of the property a) follows immediately from Propositions 1 and 3. We prove the property
b).

i) ("if-part") Let $\varkappa_0$ be a characteristic number of the
PIE (1) and $f_0\in\mathcal H$ be a solution of the PIE (1) and
$g\in \mathcal H$ be an arbitrary solution of the adjoint
homogeneous equation $f-\overline \varkappa_0S^*f=0.$ Then
$$
\langle f_0,g\rangle= \langle g_0+\varkappa_0 Sf_0,g\rangle=
\langle g_0,g\rangle+ \langle \varkappa_0 Sf_0,g\rangle= \langle
g_0,g\rangle +\langle f_0,\overline\varkappa_0 S^* g\rangle=
$$ $$
= \langle g_0,g\rangle+ \langle f_0,g\rangle.
$$
Therefore we have $\langle g_0,g\rangle=0.$

ii) ("only if"-part) Let $\varkappa_0$ be a characteristic number
of the PIE (1). Suppose that $g_0$ satisfies the condition (I),
i.e. $\langle g_0,g\rangle=0$ for every solution $g\in \mathcal H$
of the equation $f-\overline \varkappa_0S^*f=0.$

Consider the function $p(x,s,y)\in L_2(\Omega^3)$ given by the
equality
$$
p(x,s,y)= q(x,s,y)-\sum\limits_{j=1}^m \overline{f_j(s,y)}
g_j(x,y),
$$
where $f_1, f_2, \ldots, f_m$ and $g_1, g_2, \ldots, g_m$ are
orthonormal systems of the solutions of the equations (4) and
$(6'),$ respectively, in the sense of $L^0$-valued internal
product. Then for almost all $\alpha \in \Omega$ the homogeneous
Fredholm equation $\varphi-\varkappa_0 W_\alpha \varphi=0$ has in
$L_2(\Omega)$ only the trivial solution [10] where $W_\alpha$ is
an integral operator in $L_2(\Omega)$ with the kernel $p(x, s,
\alpha).$ Hence by Lemma 2 the PIE $f-\varkappa_0Wf=g_0$ has the
solution $f_0\in \mathcal H$ of the form
$$
f_0=g_0(x,y)+\varkappa_0Sf_0(x,y)- \varkappa_0 \sum\limits_{j=1}^m
\langle f_0, f_j\rangle \cdot g_j(x,y).
$$

Therefore, we obtain that
$$
\langle f_0,g_k\rangle= \langle g_0,g_k\rangle +\langle
\varkappa_0Sf_0,g_k\rangle- \sum\limits_{j=1}^{m} \langle
f_0,f_j\rangle\cdot\langle \varkappa_0g_j,g_k\rangle=
$$ $$
=\langle f_0,\overline \varkappa_0S^*g_k\rangle-\varkappa_0\langle
f_0,f_k\rangle =\langle f_0,g_k\rangle - \varkappa_0\langle
f_0,f_k\rangle,
$$
i.e. $\langle f_0,f_k\rangle=0,$ since $\varkappa_0\ne 0.$ Thus,
the solution $f_0$ of the equation $f-\varkappa_0 Wf=g_0$ has the
form $f_0=g_0+\varkappa_0Sf_0$ and hence, the function $f_0$ is
also a solution of the PIE (1) at $\varkappa=\varkappa_0.$
\end{proof}

If there exists a number $C$ such that $$|b(t)|\le C \ \mbox{ for
almost all }\ t\in\Omega, \eqno (II)$$ then the PIO $S$ is a
bounded operator on the space $L_2(\Omega^2),$ i.e $Sf\in
L_2(\Omega^2),$ $\forall f\in L_2(\Omega^2) \subset
L_{2,0}(\Omega^2)$ and $\| Sf\|_{L_2(\Omega^2)} \le C_0 \| f\|_{
L_2(\Omega^2)}$ for all $f\in L_2(\Omega^2)$ where $C_0$ is a
positive number, $b(t) = \int\int |q(x,s,t)|^2 dx ds.$

Let $k(x,s,y) \in C(\Omega^3).$ Then the subspace $L_2(\Omega^2)$
is invariant for the PIO $T_1: (T_1f) (x,y) = \int k(x,s,y) f(s,y)
ds.$ Therefore it is possible to study solvability for the PIE
\begin{equation}
f-\varkappa T_1f =g_0
\end{equation}

and it is uniquely defined by its kernel $k(x,s,y),$ where
$$
b(t)=\int\limits_\Omega \int\limits_\Omega |k(x,s,t)|^2dxds.
$$
in the space $L_2(\Omega^2)$ where $f$ is an unknown function from
$L_2(\Omega^2),$ $g_o \in L_2(\Omega^2)$ is given (known)
functionm $\varkappa \in \mathbb C$ is a parameter of the
equation.

Let $\chi_{T_1}$ be a set of characteristic numbers for the PIE
(7) (see [9]). Definition of characteristic number [9] and Theorem
3 we obtained imply
\begin{thm} Let $\varkappa_0 \in \chi_{T_1}.$ Then

a) the homogeneous equation $f-\varkappa T_1 f = \theta$ has a
non-trivial solution, moreover the set of all solutions of the
homogeneous equations is an infinite dimensional subspace of
$L_2(\Omega^2);$

b) PIE (7) is solvable if and only if a given function $g_0$
satisfies the condition
$$
\int g_0(s,y) \overline{g(s,t)} ds =0 \ \mbox{ for almost all }
t\in \Omega \eqno (III)
$$ where $g\in L_2(\Omega^2)$ is an arbitrary solution of the
adjoint homogeneous equation $f -\overline{\varkappa_0} T_1^* f =
\theta.$
\end{thm}

\end{document}